\documentclass{article}

\usepackage{amsmath,amssymb,epsf,wrapfig}
\usepackage{graphicx}
\def\qed{\hfill$\square$}

\newcommand{\R}{\mathbb{R}}

\newtheorem{theorem}{Theorem}

\newtheorem{lemma}[theorem]{Lemma}

\newtheorem{remark}[theorem]{Remark}
\newtheorem{example}[theorem]{Example}
\newtheorem{definition}[theorem]{Definition}

\begin{document}

\title{Three-dimensional braids and their descriptions}

\author{
J. Scott Carter\\
Department of Mathematics, University of South Alabama, \\
Mobile, AL 36688, USA
\footnote{E-mail address: {\tt carter@southalabama.edu} }
\and 
Seiichi Kamada\\
Department of Mathematics, Osaka City University, \\
Sumiyoshi, Osaka 558-8585, Japan
\footnote{E-mail address: {\tt skamada@sci.osaka-cu.ac.jp} }
}
\date{}

\maketitle

\begin{abstract}

The notion of a braid is generalized into two and three dimensions.  
Two-dimensional braids are described by braid monodromies or graphics called charts. 
In this paper we introduce the notion of curtains, and show that three-dimensional braids are 
described by braid monodromies or curtains. 

\end{abstract}

%
%

\section{Introduction} 

Throughout this paper, 
we work in the PL category (\cite{Hud1969, RS1972}) and assume that all manifolds are oriented and 
$m$-manifolds embedded in $(m+2)$-manifolds are locally flat.  
We denote by $D^2$ the $2$-disk and by $B^m$ the $m$-disk.  
Let $d$ be a positive integer and $X_d$ a fixed set of $d$  interior points of the $2$-disk $D^2$.  

For the product space $D^2 \times \Sigma^m$ of $D^2$ and an $m$-manifold $\Sigma_m$, 
we denote by $pr_1: D^2 \times \Sigma^m \to D^2$ the first factor projection, and by 
$pr_2: D^2 \times \Sigma^m \to \Sigma^m$  the second factor projection.

First we introduce the notion of a 3-dimensional braid.  

\begin{definition}{\rm 
(1) A {\it 3-dimensional braid}  in $D^2 \times B^3$ (or over $B^3$) of degree $d$ is a 3-manifold $M$ 
embedded in $D^2 \times B^3$   
such that (i) the restriction map $pr_2|_M : M \to B^3$ is a simple branched covering map of degree $d$ 
branched along a link in $B^3$ 
and (ii) $\partial M = M \cap \partial (D^2 \times B^3) = X_d \times \partial B^3$.  

(2) A {\it 3-dimensional braid} in $D^2 \times S^3$ (or 
over $S^3$) of degree $d$ is a 3-manifold $M$ 
embedded in $D^2 \times S^3$   
such that (i) the restriction map $pr_2|_M : M \to S^3$ is a simple branched covering map of degree $d$ branched along a link in $S^3$.   
}\end{definition} 

When we refer to a link, it may be the empty set. 
Refer to \cite{BE1979, BE1984} for simple branched coverings. 

More generally, we introduce the notion of a braided $3$-manifold as follows.  Let $\Sigma^3$ be a 3-manifold. 

\begin{definition}{\rm 
A {\it braided 3-manifold}  in $D^2 \times \Sigma^3$ (or over $\Sigma^3$) of degree $d$ 
is a 3-manifold $M$ embedded  in $D^2 \times \Sigma^3$   such that 
the restriction map $pr_2|_M : M \to \Sigma^3$ is a simple branched covering map of degree $d$ 
and $\partial M = M \cap \partial (D^2 \times \Sigma^3) \subset {\rm int} D^2 \times \partial \Sigma^3$.  
}\end{definition}

A 3-dimensional braid  in $D^2 \times B^3$  is a braided 3-manifold in $D^2 \times B^3$ 
such that $\partial M = X_d \times \partial B^3$ and the branch set is a link in $B^3$.    A 3-dimensional braid  in $D^2 \times S^3$ is a braided 3-manifold in $D^2 \times S^3$ such that the branch set is a link in $S^3$.  

Since any closed 3-manifold can be presented as a simple branched covering of $S^3$ branched along a link \cite{Hi1976, Mo1976}, our assumption that the branch set is a link is not so restrictive.  

In this paper, we study how to describe 3-dimensional braids.  We consider two methods, one is braid monodromies and the other is curtain descriptions.  The idea of the curtain description was introduced in \cite{CK2012}, and some examples were shown in \cite{CK2012, CK*}.  However, existence of a curtain for any $3$-dimensional braid was not shown.  The main purpose of this paper is to show how to construct a curtain.  

The first author was supported by the Ministry of
Education Science and Technology (MEST) and
the Korean Federation of Science and Technology
Societies (KOFST) during the initial phases of this work. The second author was supported by 
JSPS KAKENHI Grant Number 21340015. 

%
%

\section{2-dimensional braids, braid monodromies and charts} 

Before going to the case of 3-dimension in the next section, we quickly recall the notions of 2-dimensional braids, braid monodromies and charts. 
For the precise definitions and details, refer to \cite{CKS2004, Kam2002}.   The reader who is familiar with these notions may skip this section.

Let $\Sigma^2$ be a  surface.  

\begin{definition}{\rm 
A {\it braided surface}  in $D^2 \times \Sigma^2$ (or over $\Sigma^2$) of degree $d$ 
is a  surface $S$ embedded  in $D^2 \times \Sigma^2$ such that 
the restriction map $pr_2|_S : S \to \Sigma^2$ is a simple branched covering map of degree $d$ 
and $\partial S = S \cap \partial (D^2 \times \Sigma^2) \subset {\rm int} D^2  \times \partial \Sigma^2$.  
\begin{itemize}

\item[(1)] A {\it 2-dimensional braid} in $D^2 \times B^2$ (or over $B^2$)  is a braided surface in $D^2 \times B^2$ 
such that $\partial S = X_d \times \partial B^2$.    

\item[(2)] A {\it 2-dimensional braid} in $D^2 \times S^2$ (or 
over $S^2$) is a braided surface in $D^2 \times S^2$.  
\end{itemize}
}\end{definition}

\begin{definition}{\rm 
Two $2$-dimensional braids $S$ and $S'$  in $D^2 \times B^2$ are said to be {\it equivalent} if there is an ambient isotopy $\{ h_s: D^2 \times B^2 \to D^2 \times B^2 \}_{s \in [0,1]}$ such that 
\begin{itemize}

\item[(1)] $h_0 = {\rm id}$ and $h_1(S)  = S'$,   

\item[(2)] there is an ambient isotopy $\{ \underline{h}_s : B^2 \to B^2 \}_{s \in [0,1]}$ with 
$\underline{h}_s \circ pr_2 = pr_2 \circ h_s$ for each $s \in [0,1]$, and 

\item[(3)] for each $s \in [0,1]$, the restriction map of $h_s$ to $D^2 \times \partial B^2$ is the identity map. 
\end{itemize}
Moreover, if $\underline{h}_s = {\rm id}: B^2 \to B^2$ for each $s \in [0,1]$, then we say that $S$ and $S'$ are {\it isomorphic}.  
}\end{definition}

We assume that the points of $X_d$ are arranged on a straight line and identify the fundamental group $\pi_1(C_d, X_d)$
of the unordered configuration space $C_d$ of  $d$ interior points of $D^2$ with base point $X_d$ with Artin's braid group $B_d$ (cf. \cite{Kam2002}).  

Let $S$ be a $2$-dimensional braid in $D^2 \times B^2$ of degree $d$.  Take a point $q_0$ in $\partial B^2$.  
Let $\Delta(S)$ be the set of branch values of the branched covering $S \to B^2$.  

\begin{definition}{\rm 
The {\it braid monodromy} of $S$ is the homomorphism 
$$\rho_S: \pi_1(B^2 \setminus \Delta(S), q_0) \to \pi_1(C_d, X_d) = B_d$$ 
sending the homotopy class of a path $\alpha: ([0,1], \{0,1\}) \to (B^2 \setminus \Delta(S), q_0)$ to the braid presented by 
the path $([0,1], \{0,1\}) \to (C_d, X_d)$ with $t \mapsto pr_1(S \cap pr_2^{-1}(\alpha(t)))$.  
}\end{definition}

\begin{figure}[htb]
\begin{center}
\includegraphics[width=3.5in]{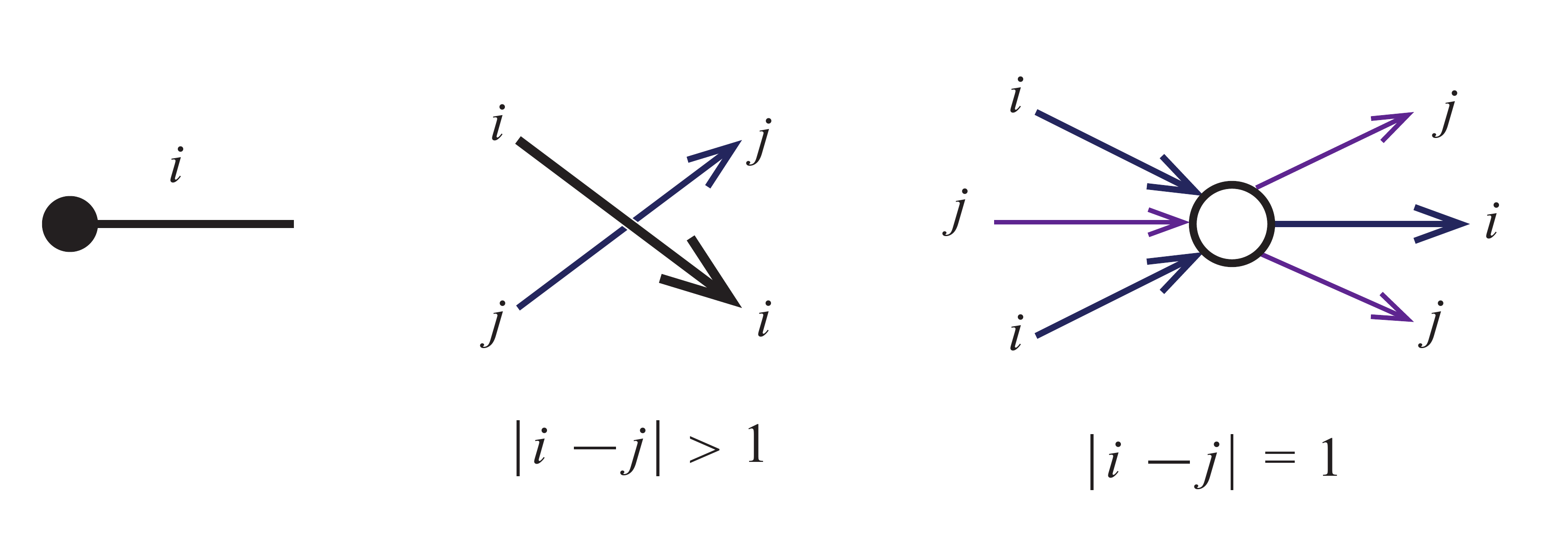}
\end{center}
\caption{Vertices of a chart}
\label{fig:schartvert_01ori}
\end{figure}

\begin{definition}{\rm 
A {\it chart} of degree $d$  is a labeled and oriented graph $\Gamma$ in $\Sigma^2$ 
such that $\Gamma \cap \partial \Sigma^2 = \emptyset$ and that 
each edge is labeled in $\{1, \dots, d-1\}$ and each vertex is as in Figure~\ref{fig:schartvert_01ori}. 
We call a vertex a {\it black vertex}, a {\it crossing} or a {\it white vertex} if the valency of the vertex is $1$, $4$ or $6$, respectively.  The arrow at a black vertex in this figure is suppressed since it may either be incoming or outgoing. 

When $\partial \Sigma^2 \neq \emptyset$, a {\it chart with external boundary} is a chart for which we allow the case $\Gamma \cap \partial \Sigma^2 \neq \emptyset$ such that the intersection $\Gamma \cap \partial \Sigma^2$ consists of degree-$1$ vertices of $\Gamma$. We call the degree-$1$ vertices on $\partial \Sigma^2$ the {\it boundary vertices} or {\it external boundary vertices} of $\Gamma$.  (Note that degree-$1$ vertices in ${\rm int} \Sigma^2$ are called black vertices.) 

}\end{definition}

Let $\sigma_i$ $(i=1, \dots, d-1)$ be the standard generators of the braid group $B_d$.  Then $B_d$ has a group presentation 
$$B_d = \left\langle \sigma_1, \dots, \sigma_{d-1}  \, 
\begin{array}{|ll}
\sigma_i \sigma_j \sigma_i = \sigma_j \sigma_i \sigma_j  \quad & (|i-j| = 1) \\ 
\sigma_i \sigma_j = \sigma_j \sigma_i &   (|i-j| > 1) 
\end{array}
\right\rangle. $$

Let $\Gamma$ be a chart in $B^2$ of degree $d$. 
 Take a point $q_0$ in $\partial B^2$.  
Let $\Delta(\Gamma)$ be the set of black vertices of $\Gamma$.   

For a path $\alpha: ([0,1], \{0,1\}) \to (B^2 \setminus \Delta(\Gamma), q_0)$, up to homotopy, we assume that $\alpha$ intersects with $\Gamma$ in general position.  We associate with a letter $\sigma_{i_k}^{\epsilon_k}$ the $k$th intersection of $\alpha$ with $\Gamma$, where $i_k$ is the label on the edge of $\Gamma$ and $\epsilon_k$ is $+1$ or $-1$ determined by the orientations of $\alpha$ and the edge.  Reading these letters along $\alpha$ we have a word on the standard generators of $B_d$.   The word is called the {\it intersection word} of $\alpha$ with respect to $\Gamma$.  

\begin{definition}{\rm 
The {\it braid monodromy} of $\Gamma$ is the homomorphism 
$$\rho_\Gamma: \pi_1(B^2 \setminus \Delta(\Gamma), q_0) \to  B_d$$ 
sending the homotopy class of a path $\alpha: ([0,1], \{0,1\}) \to (B^2 \setminus \Delta(\Gamma), q_0)$ to the braid presented by 
the intersection word of $\alpha$ with respect to $\Gamma$.   
}\end{definition}

\begin{figure}[htb]
\begin{center}
\includegraphics[width=4in]{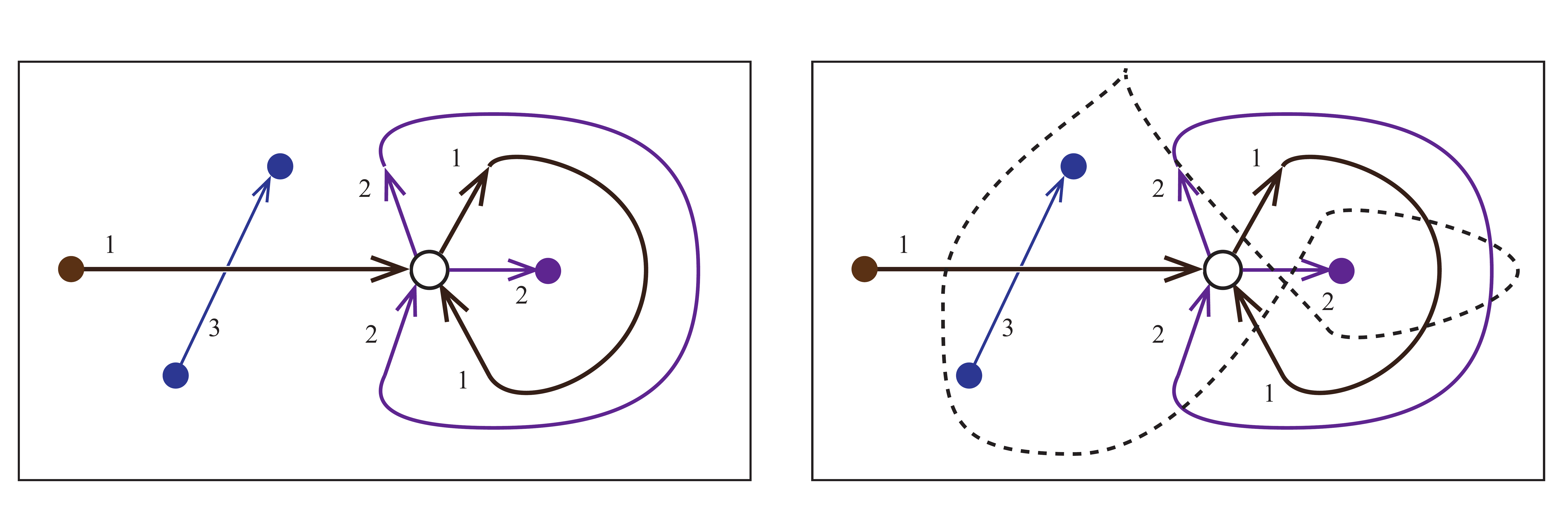}
\end{center}
\caption{A chart}
\label{fig:scharteg_03bCS3}
\end{figure}

Note that, by the correspondence $i  \leftrightarrow  \sigma_i  \in B_d$, the labels of a chart are assumed to present the standard generators in $B_d$.   The left of Figure~\ref{fig:scharteg_03bCS3} is an example of a chart of degree~$4$. 
For the path $\alpha$ illustrated in Figure~\ref{fig:scharteg_03bCS3} as a dotted line, 
the intersection word is $\sigma_1^{-1} \sigma_2^{-1} \sigma_1^{-1} \sigma_2 \sigma_2 \sigma_1 \sigma_2$, 
and $\rho_\Gamma([\alpha]) = \sigma_1 \in B_4$.  

\begin{theorem}[\cite{Kam1992}, cf. \cite{Kam2002}]\label{2braid-chart}
For  any 2-dimensional braid $S$ in $D^2 \times B^2$, 
there exists a chart $\Gamma$ with $\rho_S = \rho_\Gamma$. 
Conversely, for any chart $\Gamma$, there exists a 2-dimensional braid $S$ in $D^2 \times B^2$ with 
$\rho_S = \rho_\Gamma$.  
\end{theorem} 

In the situation of this theorem, we call $\Gamma$ a {\it chart description} of $S$. 

\begin{figure}[htb]
\begin{center}
\includegraphics[width=4.5in]{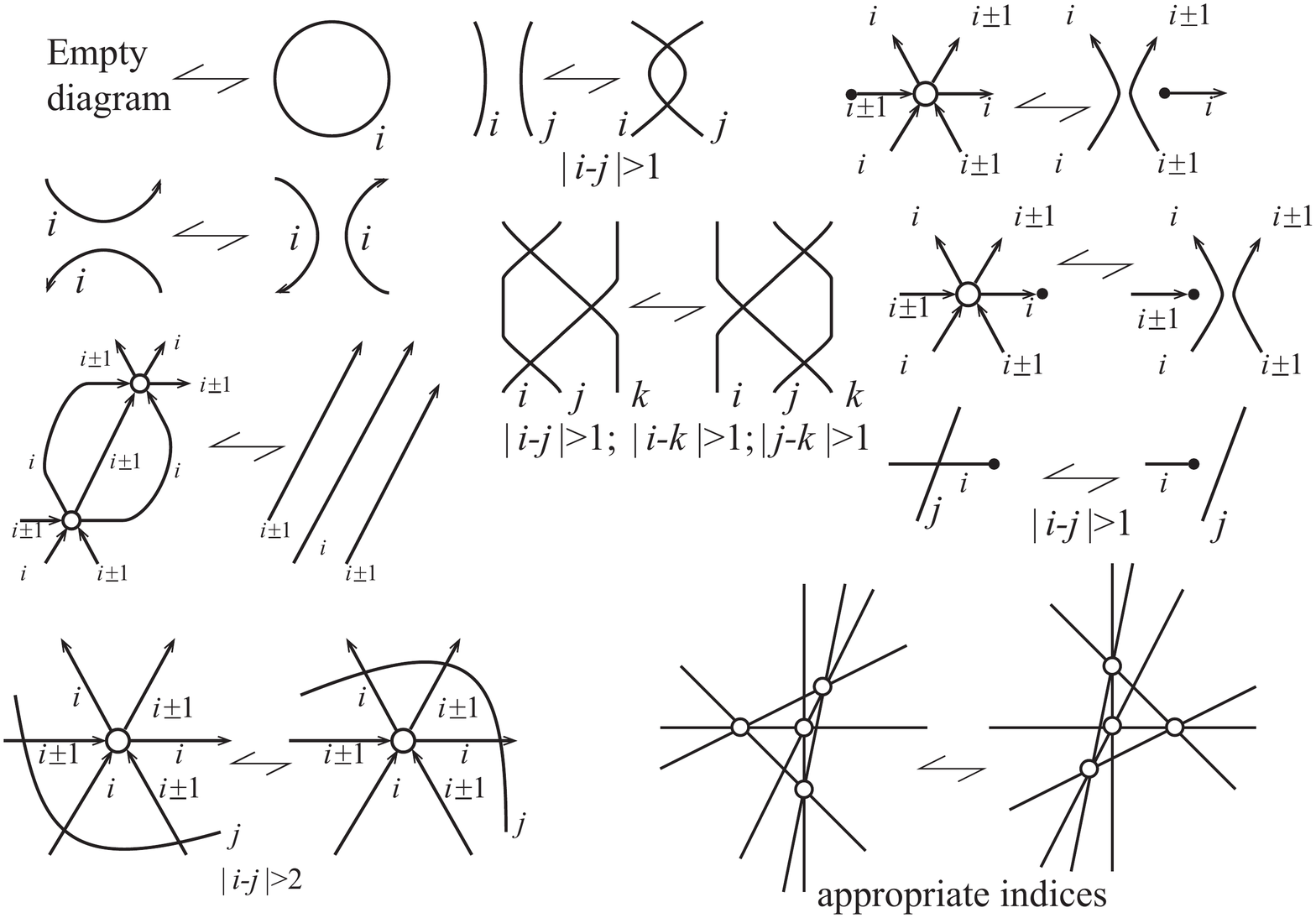}
\end{center}
\caption{Chart moves}
\label{fig:schartmovesCS3}
\end{figure}

The local moves on charts illustrated in Figure~\ref{fig:schartmovesCS3} are called {\it chart moves}.  

\begin{theorem}[\cite{Kam1992, Kam1996}, cf. \cite{CS1998, Kam2002, Ta2007}] 
Let $S$ and $S'$ be 2-dimensional braids in $D^2 \times B^2$, and let $\Gamma$ and $\Gamma'$ be chart descriptions of them, respectively.  
\begin{itemize}
\vspace{-0.2cm}
\item[(1)] $S$ and $S'$ are isomorphic if and only if $\Delta(\Gamma)= \Delta(\Gamma')$ and 
$\Gamma$ is related to $\Gamma'$ by a finite sequence of ambient isotopies of $B^2$ and chart moves keeping 
$\Delta(\Gamma)$ fixed.  
\vspace{-0.2cm}
\item[(2)]  
$S$ and $S'$ are equivalent if and only if  
$\Gamma$ is related to $\Gamma'$ by a finite sequence of ambient isotopies of $B^2$ and chart moves.  
\end{itemize}
\end{theorem} 

%
%

\section{3-dimensional braids } 

First we introduce the notion of equivalence on 3-dimensional braids.  

\begin{definition}\label{defn:3baidequiv}{\rm 
Two $3$-dimensional braids $M$ and $M'$  in $D^2 \times B^3$ are said to be {\it equivalent} if there is an ambient isotopy $\{ h_s: D^2 \times B^3 \to D^2 \times B^3 \}_{s \in [0,1]}$ such that 
\begin{itemize}

\item[(1)] $h_0 = {\rm id}$ and $h_1(M)  = M'$,   

\item[(2)] there is an ambient isotopy $\{ \underline{h}_s : B^3 \to B^3 \}_{s \in [0,1]}$ with 
$\underline{h}_s \circ pr_2 = pr_2 \circ h_s$ for each $s \in [0,1]$, and 

\item[(3)] for each $s \in [0,1]$, the restriction map of $h_s$ to $D^2 \times \partial B^3$ is the identity map. 

\end{itemize}
Moreover, if $\underline{h}_s = {\rm id}: B^3 \to B^3$ for each $s \in [0,1]$, then we say that $M$ and $M'$ are {\it isomorphic}.  
}\end{definition}

Let $M$ and $M'$ be $3$-dimensional braids in $D^2 \times B^3$ of degree $d$, and let 
$L= \Delta(M)$ and $L'= \Delta(M')$ be the set of branch values of the branched coverings 
$M \to B^3$ and $M' \to B^3$, respectively.  If $M$ and $M'$ are equivalent, then $L$ is ambient isotopic to $L'$ in $B^3$. 
If $M$ and $M'$ are isomorphic, then $L=L'$.  

The notion of a braid monodromy is also defined for 3-dimensional braids.  

Let $M$ be a $3$-dimensional braid in $D^2 \times B^3$ of degree $d$.  Take a point $q_0$ in $\partial B^3$.  
Let $\Delta(M)$ be the set of branch values of the branched covering $M \to B^3$.

\begin{definition}{\rm 
The {\it braid monodromy} of $M$ is the homomorphism 
$$\rho_M: \pi_1(B^3 \setminus \Delta(M), q_0) \to \pi_1(C_d, X_d) = B_d$$ 
sending the homotopy class of a path $\alpha: ([0,1], \{0,1\}) \to (B^3 \setminus \Delta(M), q_0)$ to the braid presented by 
the path $([0,1], \{0,1\}) \to (C_d, X_d)$ with $t \mapsto pr_1(M \cap pr_2^{-1}(\alpha(t)))$.  
}\end{definition}

\begin{figure}[htb]
\begin{center}
\includegraphics[width=1.5in]{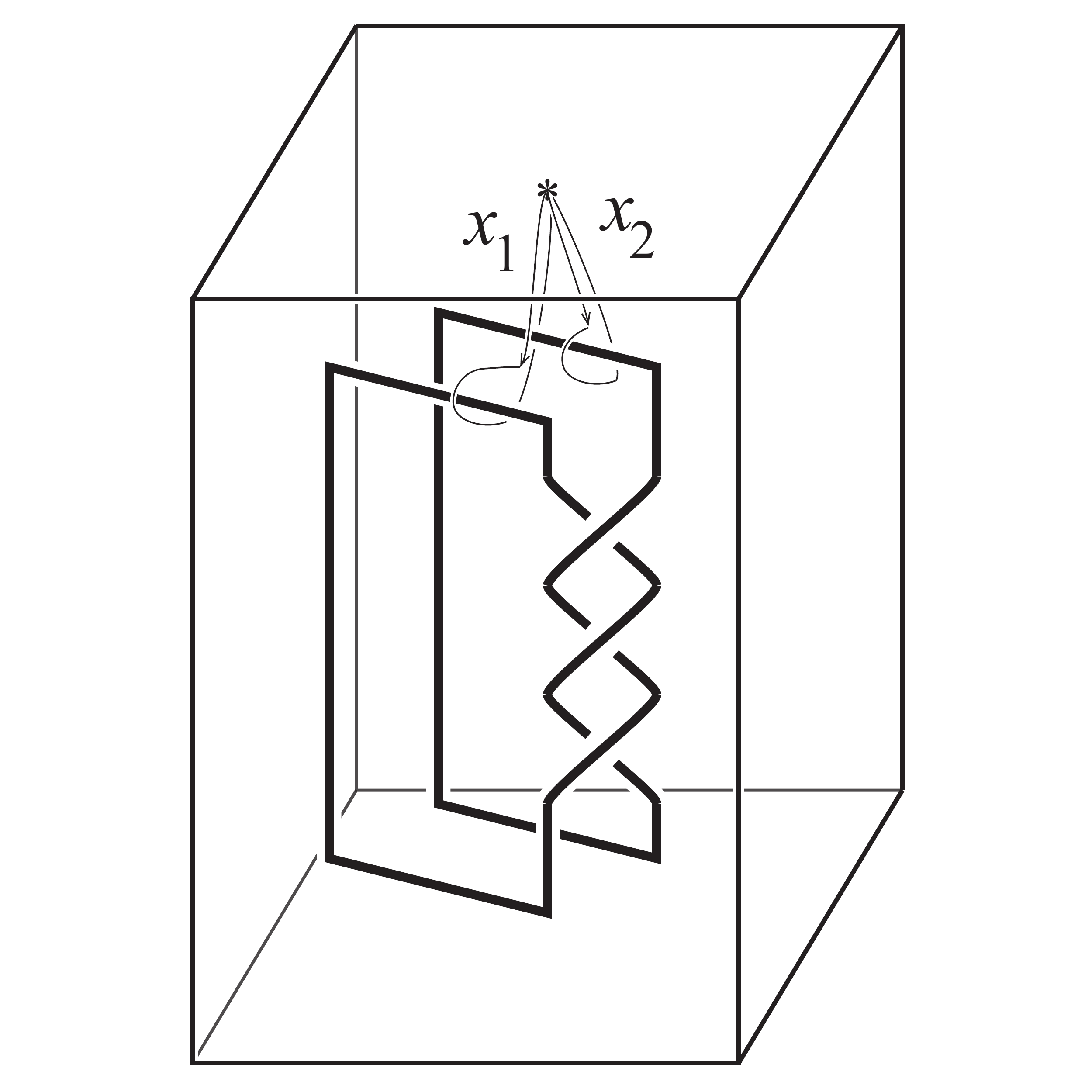}
\end{center}
\caption{A branch set $L= \Delta(M)$}
\label{fig:s3_1}
\end{figure}

\begin{example}{\rm 
Let $L$ be a trefoil in $B^3$ as in Figure~\ref{fig:s3_1}.  
There is a homomorphism $\rho: \pi_1(B^3 \setminus L, q_0) \to B_3$ with  
$\rho(x_1)= \sigma_1$ and $\rho(x_2) = \sigma_2$.  
There is a $3$-dimensional braid $M$ in $D^2 \times B^3$ of degree $3$ with $\rho_M= \rho$.  We will show how to construct such a $3$-dimensional braid $M$ later.  
}\end{example}

Two braid monodromies  
$\rho: \pi_1(B^3 \setminus \Delta, q_0) \to B_d$ and  
$\rho': \pi_1(B^3 \setminus \Delta' q_0) \to B_d$ are said to be 
{\it equivalent} 
if there is a 
homeomorphism $h: (B^3, \Delta) \to (B^3, \Delta')$ rel $\partial B^3$ 
such that $\rho= \rho' \circ h_\ast$, where 
$h_\ast : \pi_1(B^3 \setminus \Delta, q_0) \to \pi_1(B^3 \setminus \Delta', q_0)$ is the isomorphism induced from $h$.  

\begin{lemma}
Let  $M$ and $M'$ be $3$-dimensional braids in $D^2 \times B^3$ of degree $d$.  
\begin{itemize}
\item[(1)] If $M$ and $M'$ are equivalent, then $\rho_M$ and $\rho_{M'}$ are equivalent.   
\item[(2)] If $M$ and $M'$ are isomorphic, then $\rho_M = \rho_{M'}$.  
\end{itemize}
\end{lemma}

{\it Proof.} Let $\{ h_s : D^2 \times B^3 \to D^2 \times B^3 \}_{s \in [0,1]}$ and 
 $\{ \underline{h}_s : B^3 \to B^3 \}_{s \in [0,1]}$ be ambient isotopies as in Definition~\ref{defn:3baidequiv}.  
Put $ h = \underline{h}_1: B^3 \to B^3$, which is a homeomorphism  rel $\partial B^3$.  
Since $h_1(M)=M'$ and 
$h \circ pr_2|_M = pr_2|_{M'} \circ h_1|_M$, we have 
$h(\Delta(M))= \Delta(M')$ and $\rho= \rho' \circ h_\ast$.  \qed 

\begin{lemma}\label{lem:3baridmonodromy-character}
Let $\rho_M : \pi_1(B^3 \setminus \Delta(M), q_0) \to B_d$ be the braid monodromy of a $3$-dimensional braid $M$ of degree $d$.  Then, for any meridian element $x$ of $\pi_1(B^3 \setminus \Delta(M), q_0)$, $\rho_M(x)$ is a conjugate of a standard generator of $\sigma_k$ or $\sigma_k^{-1}$.  
\end{lemma}

{\it Proof.} 
Let $D$ be a meridian disk of the link $\Delta(M)$ in $B^3$ such that $x$ is represented by $\beta^{-1} \cdot \partial D \cdot \beta$, where $\beta$ is a path in $B^3 \setminus \Delta(M)$ connecting $q_0$ and a point of $\partial D$.  The preimage $(pr_2|_M)^{-1}(D)$ is the disjoint union of $d-1$ $2$-disks in $M$ and one of them, say $D'$, is mapped onto $D$ as a $2$-fold branched covering.  The boundary $\partial D'$ is a $(2,1)$-torus knot in the solid torus $D^2 \times \partial D = pr_2^{-1}(\partial D)$, and $\rho_M(x)$ is a conjugate of a standard generator of $\sigma_k$ or $\sigma_k^{-1}$ (cf. Lemma~16.12 of \cite{Kam2002}).  \qed

\begin{figure}[htb]
\begin{center}
\includegraphics[width=3in]{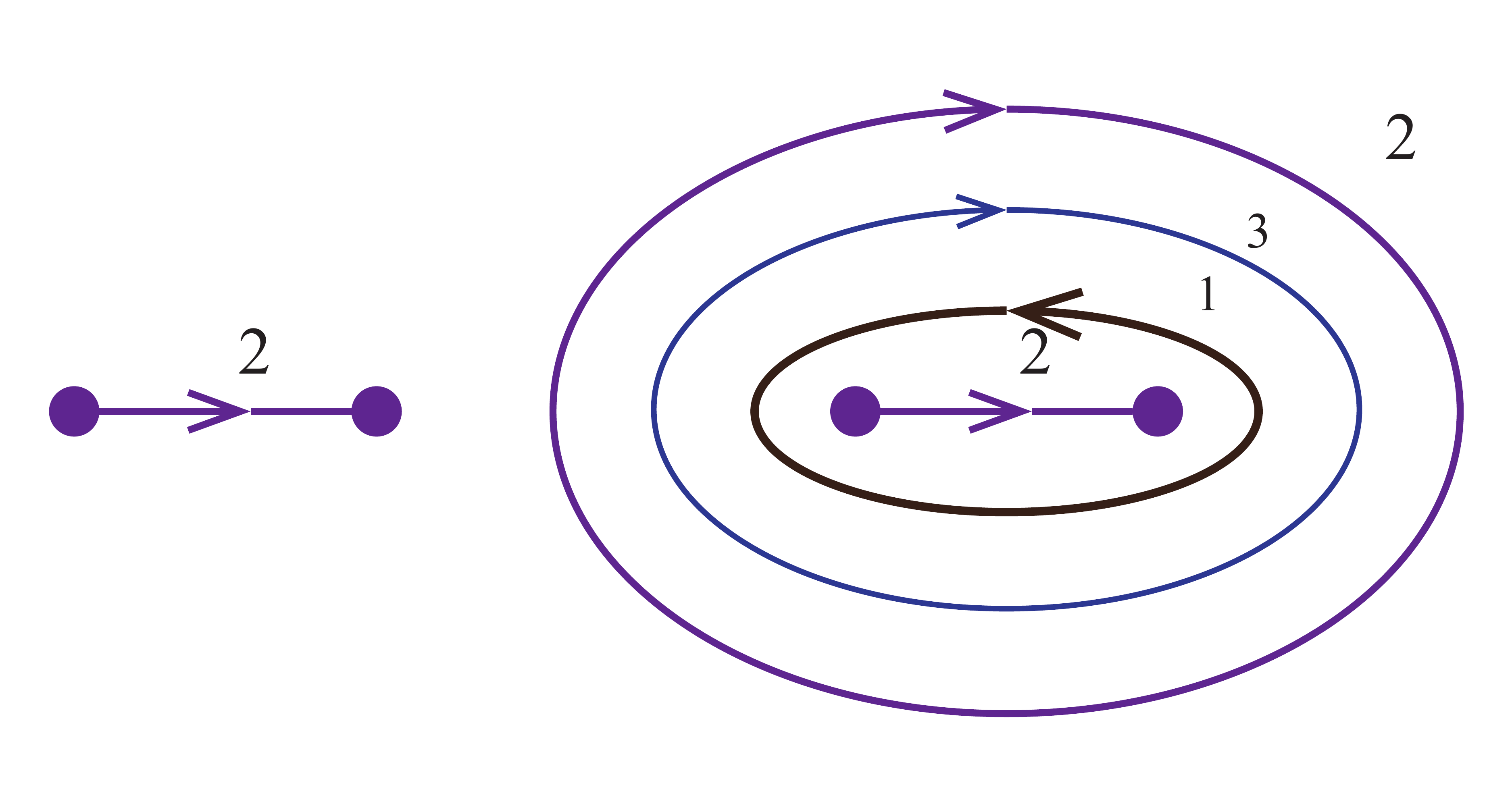}
\end{center}
\caption{A free edge and an oval nest}
\label{fig:sfedge}
\end{figure}

For a chart in $B^2$, an edge both of whose endpoints are black vertices is called a {\it free edge}. 
An {\it oval nest} is a free edge together with some concentric simple loops (\cite{Kam2002}).  
See Figure~\ref{fig:sfedge}; the left is a free edge and the right is an oval nest.  

We identify $B^3$ with $B^2 \times [0,1]$.  For a subset $X$ in $B^3 = B^2 \times [0,1]$, the {\it motion picture} of $X$ means the one-parameter family $\{ X_t \}_{t \in [0,1]}$ with $X \cap B^2 \times \{t\} = X_t \times \{t\}$ for $t \in [0,1]$.  
When $X$ is an oriented 2-simplex, $X_t$ is an oriented interval, a point, or the empty set in $B^2$.  

\begin{definition}\label{defn:curtain}{\rm 
A {\it curtain} in $B^3$ of degree $d$  is a $2$-complex $C$ in $B^3$ 
such that $C \cap \partial B^3 = \emptyset$ and that 
each face is labeled in $\{1, \dots, d-1\}$ and oriented, and after a slight perturbation,  the motion picture $\{ C_t \}_{t \in [0, 1]}$ of $C$ satisfies the following conditions: 
\begin{itemize}

\item[(1)] The slices $C_t$ are charts of degree $d$ for all but a finite number of $t$. Let $t_1, \dots, t_k$ be the exceptional values. And put $t_0=0, t_{k+1}=1$.  

\item[(2)] For each interval $(t_i, t_i+1)$ $(i=0, \dots, k)$, $C_t$ is deformed by an ambient isotopy of $B^2$ rel $\partial B^2$. 

\item[(3)]  For each exceptional value $t_i$ $(i=1, \dots, k)$, a chart move, insertion of some free edges or  deletion of some free edges occurs. 
\end{itemize}

A {\it curtain with external boundary} is a curtain for which we allow the case $C \cap \partial B^3 \neq \emptyset$ such that the slices $C_t$ are charts with external boundary.  ($C_0$ and $C_1$ may be nonempty charts with external boundary.) 
}\end{definition}

Let $C$ be a curtain (possibly with external boundary) and let $\{ C_t \}_{t \in [0,1]}$ the motion picture satisfying the conditions in Definition~\ref{defn:curtain}.  Let $E \subset [0,1]$ be the set of exceptional values $\{ t_1, \dots, t_k\}$. 
We call the subset $\partial^{\rm int}(C):= (\cup_{t \in [0,1]} \partial^{\rm int}(C_t)) \cup (\cup_{t \in E} \{ \mbox{ free edges in $C_t$ that are inserted or deleted at $t$} \})$ 
the {\it internal boundary} of $C$, and the subset 
$\partial^{\rm ext}(C) := \cup_{t \in [0,1]} \partial^{\rm ext}(C_t)$ the {\it external boundary} of $C$, where $\partial^{\rm int}(C_t)$ is the set of black vertices of $C_t$ and $\partial^{\rm ext}(C_t)$ is the set of boundary vertices of $C_t$.  
Note that the internal boundary is a link in $B^3$ if $C_0$ and $C_1$ are the empty charts. We also denote the  internal boundary $\partial^{\rm int}(C)$ by $\Delta(C)$.

Let $C$ be a curtain in $B^3$  of degree $d$. 
 Take a point $q_0$ in $\partial B^3$.  
Let $\Delta(C)= \partial^{\rm int}(C)$ be the internal boundary of $C$. 

For a path $\alpha: ([0,1], \{0,1\}) \to (B^3 \setminus \Delta(C), q_0)$, up to homotopy, we assume that $\alpha$ intersects with $C$ in general position, namely, each intersection is a transverse intersection at an interior point of a face of $C$. 
We associate with a letter $\sigma_{i_k}^{\epsilon_k}$ the $k$th intersection of $\alpha$ with $C$, where $i_k$ is the label on the face of $C$ and $\epsilon_k$ is $+1$ or $-1$ determined by the orientations of $\alpha$ and the face.  Reading these letters along $\alpha$ we have a word on the standard generators of $B_d$.   The word is called the {\it intersection word} of $\alpha$ with respect to $C$.  

\begin{definition}{\rm 
The {\it braid monodromy} of $C$ is the homomorphism 
$$\rho_C: \pi_1(B^3 \setminus \Delta(C), q_0) \to  B_d$$ 
sending the homotopy class of a path $\alpha: ([0,1], \{0,1\}) \to (B^3 \setminus \Delta(C), q_0)$ to the braid presented by 
the intersection word of $\alpha$ with respect to $C$.   
}\end{definition}

\begin{theorem}\label{thm:monodromy-curtain}
Let $L$ be a link in $B^3$. 
For any homomorphism $\rho: \pi_1(B^3 \setminus L, q_0) \to B_d$ sending each meridian element to a conjugate of $\sigma_k$ or $\sigma_k^{-1}$, there exists a curtain $C$ of degree $d$ such that $\Delta(C)= L$ and 
$\rho_C = \rho$. 
\end{theorem}

We will prove this theorem in Section~\ref{sect:construct-curtain}.  

Combining this theorem and Lemma~\ref{lem:3baridmonodromy-character}, we have the following. 

\begin{theorem}\label{thm:3braid-curtain}
For any $3$-dimensional braid $M$ in $D^2 \times B^3$ of degree $d$, there exists a curtain $C$ of degree $d$ such that $\Delta(C)= \Delta(M)$ and 
$\rho_C = \rho_M$. 
\end{theorem}

In the situation of this theorem, we call $C$ a {\it curtain description} or a {\it chart description} of $M$.

%
%

\section{Construction of a curtain}\label{sect:construct-curtain}

In this section we prove Theorem~\ref{thm:monodromy-curtain}.  
For simplicity of the argument, we identify $B^3$ with $B^2 \times [-2, 2]$ (not $B^2 \times [0,1]$) and assume 
$B^2 = [-1,1] \times [0,1] \subset \R^2$.  
Put $B^-:= [-1,0] \times [0,1]$ and $B^+:= [0,1]\times [0,1]$ so that $B^2 = B^- \cup B^+$.  
Take the base point $q_0 = (q_0^\ast, 1/2) \in B^3=B^2 \times [-2,2]$ where $q_0^\ast =(0,1) \in \partial B^2 \subset \R^2$.   

Let $L$ be a link in $B^3$ and $\rho: \pi_1(B^3 \setminus L, q_0) \to B_d$ a homomorphism sending each meridian element to a conjugate of $\sigma_k$ or $\sigma_k^{-1}$.   

Assume that $L$ can be presented as a closed braid of degree $n$ for some positive integer $n$.  (By the Alexander theorem, there exists such an $n$.  Here we do not assume that $n$ is the minimum among such integers.)

\begin{figure}[htb]
\begin{center}
\includegraphics[width=3in]{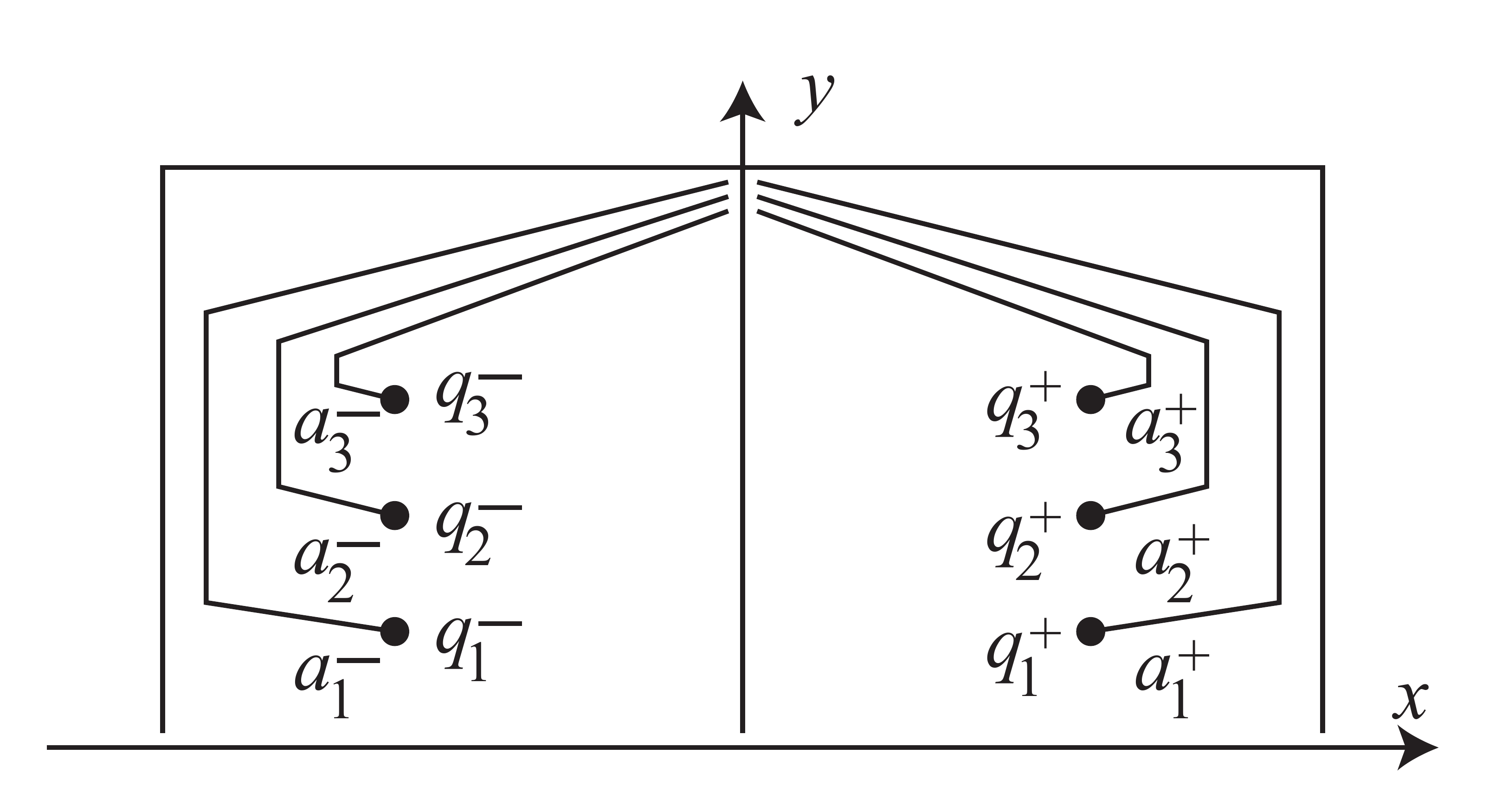}
\end{center}
\caption{Points $q_i^-$ and $q_i^+$ and arcs  $a_i^-$ and $a_i^+$}
\label{fig:shurwitzCS3}
\end{figure}

Take $n$ real numbers $q_1, \dots, q_n$ with $0 < q_1 <q_2 < \cdots < q_n <1$, and let 
$q_i^- = (-1/2, q_i) \in B^-$ and $q_i^+ = (1/2, q_i) \in B^+$ for $i= 1, \dots, n$.  
See Figure~\ref{fig:shurwitzCS3}.  (The arcs $a_i^-$ and $a_i^+$ in the figure will be used later.) 
Let $A_i$ $(i=1, \dots, n)$ be the straight segment $| q_i^- q_i^+ |$.  

By an ambient isotopy of $B^3 =B^2 \times [-2,2]$, deform the link $L$ to a link $L'$ in a braid form of degree $n$ 
satisfying the following. 
\begin{itemize}
\item[(1)] $L' \cap B^2 \times \{t \} = \emptyset$ for $t$ with $|t| > 1$.  
\item[(2)] $L' \cap B^2 \times \{t \} = \cup_{i=1}^n A_i \times \{t \}$ for $t$ with $|t| =1$.  
\item[(3)] $L' \cap B^2 \times \{t \} = \cup_{i=1}^n \{ q_i^+, q_i^-\} \times \{t\}$ for $t$ with $1/2 < |t| <1$.  
\item[(4)] $L' \cap B^+ \times [0,1/2]$ is an $n$-braid. 
\item[(5)] $L' \cap B^+ \times [-1/2,0]$ is the trivial $n$-braid $\cup_{i=1}^n \{ q_i^+\} \times [-1/2,0]$. 
\item[(6)] $L' \cap B^- \times [-1/2,1/2]$ is the trivial $n$-braid $\cup_{i=1}^n \{ q_i^-\} \times [-1/2,1/2]$. 
\end{itemize}
Let $f: B^3 \to B^3$ be a homeomorphism with $f(L')=L$ and $f |_{\partial B^3}= {\rm id}$, and let 
$\rho' = \rho \circ f_\ast: \pi_1(B^3 \setminus L', q_0) \to B_d$, where $f_\ast:  \pi_1(B^3 \setminus L', q_0) 
\to \pi_1(B^3 \setminus L, q_0)$ is the isomorphism induced from $f$.  If there exists a curtain $C'$ with $\rho_{C'}= \rho'$, then $C:= f(C)$ is a curtain with $\rho_{C} = \rho$.  

Therefore, without loss of generality, we may assume that $L$ is in a braid form as $L'$ stated above.  
Now we will construct a desired curtain $C$. 

Let $\{ L_t \}_{t \in [-2,2]}$ be the motion picture of $L \subset B^2 \times [-2,2]$.  
At the level of $t=1/2$, $L_{1/2} =  \cup_{i=1}^n \{ q_i^+, q_i^-\}$ and $\pi_1(B^2 \setminus L_{1/2}, q_0^\ast)$ is a free group generated by meridian elements.  Consider a Hurwitz arc system $(a_1^-, a_2^-, \dots, a_n^-, a_n^+, \dots, a_2^+, a_1^+)$ for $L_{1/2}$, where 
$a_i^-$ and $a_i^+$ are arcs as in Figure~\ref{fig:shurwitzCS3}.  Let $(x_1, x_2, \dots, x_n, y_n, \dots, y_2, y_1)$ be the corresponding generating system of $\pi_1(B^2 \setminus L_{1/2}, q_0^\ast)$.  Let $\iota: 
\pi_1(B^2 \setminus L_{1/2}, q_0^\ast) \to \pi_1(B^3 \setminus L, q_0)$ be the homomorphism induced from the inclusion map. 
Since $L \cap B^2 \times [1/2, 2]$ is a trivial $n$ tangle, we have  $\iota (x_i) = \iota (y_i)^{-1}$ for each $i=1, \dots, n$.  
By the assumption of $\rho$, each $\rho(\iota (x_i))$ is a conjugate of $\sigma_k$ or $\sigma_k^{-1}$.  
Therefore, there exists a chart $\Gamma$ in $B^2$ with 
$\Delta(\Gamma) = L_{1/2}$ and $\rho_\Gamma = \rho \circ \iota: \pi_1(B^2 \setminus L_{1/2}, q_0^\ast) \to B_d$.  
Moreover, by the argument in \cite{Kam2002}, this chart can be taken as a disjoint 
union of oval nests $U_1$, \dots, $U_n$ such that the free edge of the oval nest $U_i$ is $A_i$.  
(Such a chart is called a {\it ribbon chart} (\cite{Kam1992, Kam2002}).) 
We assume that each oval nest $U_i$ is symmetric with respect to the $y$-axis. 

We define $C \cap B^2 \times [1/2, 2]$ and $C \cap B^2 \times [-2, -1/2]$ as follows. 
Let $\{ C_t \}_{t \in [-2,2]}$ denote the motion picture for $C$. 

For $t$ with $1/2 \leq |t| \leq 1$, define $C_t = \cup_{i=1}^n U_i$.  
For $t$ with $1 < |t| \leq 3/2$, define $C_t = \cup_{i=1}^n U_i \setminus A_i$.  
Since $C_{3/2} =  \cup_{i=1}^n U_i \setminus A_i$ consists of simple loops, we can remove the loops one by one by  chart moves to obtain the empty chart.  By this sequence, we construct $\{ C_t \}_{t \in [3/2, 2]}$ with $C_2= \emptyset$.  
Similarly, we construct   $\{ C_t \}_{t \in [-2, -3/2]}$ with $C_{-2}= \emptyset$.  

Now we construct $C \cap B^2 \times [0, 1/2]$.  
From the configuration of $L$, there is an ambient isotopy of $B^2$ rel $(\partial B^+) \cup B^-$ 
such that the trace of the $2n$ point $\cup_{i=1}^n \{ q_i^+, q_i^-\}$ forms the braid $L \cap B^2 \times [0, 1/2]$.  
Using this isotopy, we can construct $C \cap B^2 \times [0, 1/2]$ such that the motion picture 
$\{ C_t \}_{ t \in [0, 1/2]}$ is a deformation of $C_{1/2}$ by 
an ambient isotopy of $B^2$ rel $(\partial B^+) \cup B^-$.  
Then $C_0$ is a chart which is ambient isotopic to $C_{1/2}$ in $B^2$ rel $(\partial B^+) \cup B^-$. 
Note that $\Delta(C_0) = L_0= L_{1/2} = \cup_{i=1}^n \{ q_i^+, q_i^-\}$.   

The braid monodromy $\rho_{C_0}$ is equal to $\rho_{C_{-1/2}}$.  Therefore the chart $C_0$ is related to the chart $C_{-1/2}$ by a finite sequence of ambient isotopies of $B^2$ and chart moves keeping $L_0$ fixed.  Using this sequence, we can construct a motion picture $\{ C_t \}_{ t \in [-1/2, 0]}$ connecting $C_0$ and $C_{-1/2}$.  

Now we have constructed a curtain $C$ in $B^3= B^2 \times [-2,2]$ with $\Delta(C)=L$. 
By the construction, $\rho_C = \rho$. \qed 

\begin{figure}[htb]
\begin{center}
\includegraphics[width=4in]{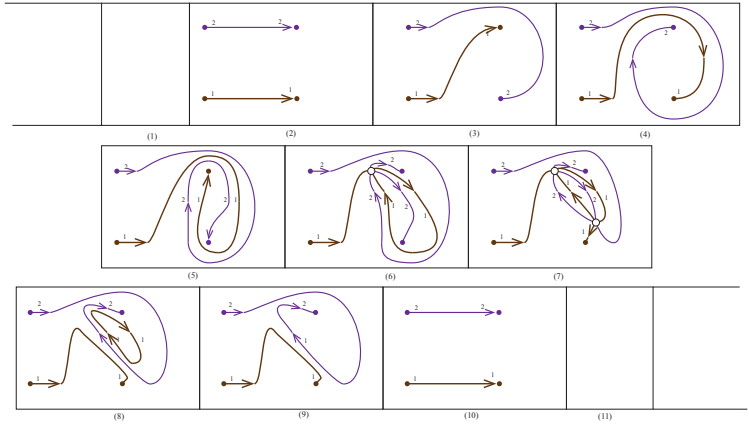}
\end{center}
\caption{A curtain}
\label{fig:scurtain}
\end{figure}

\begin{example}{\rm 
Figure~\ref{fig:scurtain} is an example of a curtain constructed by the argument in the proof of Theorem~\ref{thm:monodromy-curtain}.  Let $L$ be a trefoil in $B^3 = B^2 \times [-2,2]$ as in Figure~\ref{fig:s3_1}.  
Consider a homomorphism $\rho: \pi_1(B^3 \setminus L, q_0) \to B_3$ with  
$\rho(x_1)= \sigma_1$ and $\rho(x_2) = \sigma_2$.  
Then the curtain $M$  illustrated in Figure~\ref{fig:scurtain} satisfies that $\Delta(C)= L$ and 
$\rho_C = \rho$.  

Here is an explanation on the motion picture $\{ C_t \}_{ t \in [-2,2]}$.  
For $t \in (1, 2]$, $C_t$ is the empty chart, see (1) of Figure~\ref{fig:scurtain}.  For $t \in [1/2, 1]$, $C_t$ is a chart consisting of two free edges $A_1$ and $A_2$ as in (2) of the figure.  
For $t \in [0, 1/2]$, the motion picture $\{C_t\}$ is constructed by using an ambient isotopy of $B^2$.  
See (2)--(5).  $C_0$ is depicted in (5).  
We prepare the chart $C_{-1/2}$ as the same as $C_{1/2}$. This is (10) of the figure.  Since the chart $C_0$ and $C_{-1/2}$ describe the same braid monodromy, there is a finite sequence of chart moves changing $C_0$ to $C_{-1/2}$.  This process is depicted from (5) to (10).  For $t \in [-1, -1/2]$, $C_t$ is the same with $C_{1/2}$.  For $t \in [-2, -1)$, $C_t$ is the empty chart.  
}\end{example}

\begin{figure}[htb]
\begin{center}
\includegraphics[width=4in]{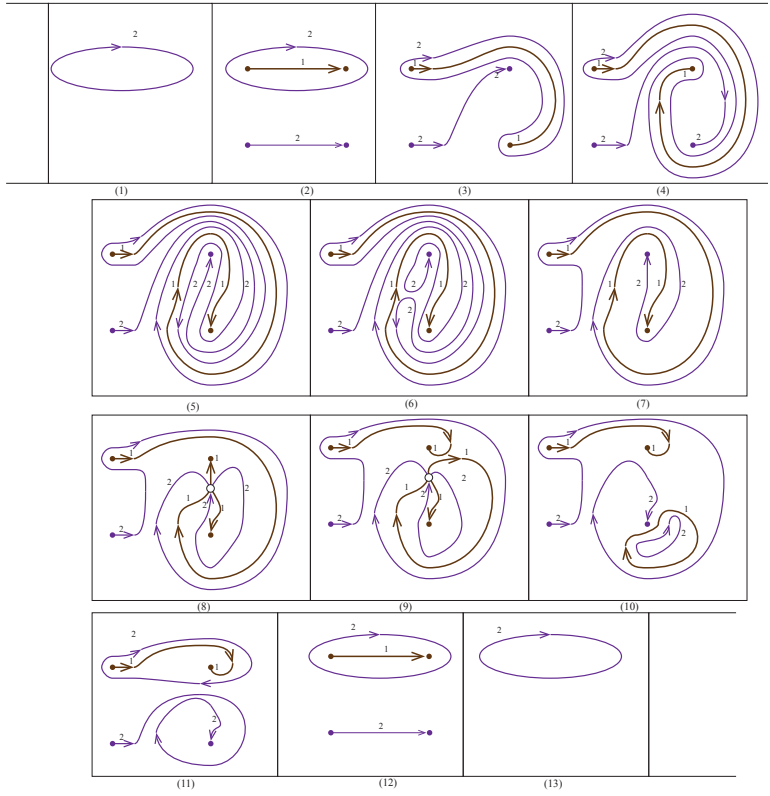}
\end{center}
\caption{A curtain}
\label{fig:scurtain2}
\end{figure}

\begin{example}{\rm 
Here is another example. 
Let $L$ be a trefoil in $B^3 = B^2 \times [-2,2]$ as in Figure~\ref{fig:s3_1}.  
Consider a homomorphism $\rho: \pi_1(B^3 \setminus L, q_0) \to B_3$ with  
$\rho(x_1)= \sigma_2$ and $\rho(x_2) = \sigma_2^{-1} \sigma_1 \sigma_2$.  
Then the curtain $M$  illustrated in Figure~\ref{fig:scurtain2} satisfies that $\Delta(C)= L$ and 
$\rho_C = \rho$.  

Here is an explanation on the motion picture $\{ C_t \}_{ t \in [-2,2]}$.  
For $t \in (1, 3/2]$, $C_t$ is a chart consisting of a simple loop as in (1) of the figure. 
The motion picture of $C \cap B^2 \times [3/2, 2]$ has a chart move removing the loop of 
$C_{3/2}$.  
For $t \in [1/2, 1]$, $C_t$ is a chart consisting of a free edge and an oval nest as in (2) of the figure.  
For $t \in [0, 1/2]$, the motion picture $\{C_t\}$ is a deformation of $C_{1/2}$ by 
an ambient isotopy of $B^2$.  $C_0$ is depicted in (5) of the figure.  
We prepare the chart $C_{-1/2}$ as the same as $C_{1/2}$. This is (12) of the figure. 
Since the chart $C_0$ and $C_{-1/2}$ describe the same braid monodromy, there is a finite sequence of chart moves changing $C_0$ to $C_{-1/2}$.  This process is depicted from (5) to (12) in the figure.  For $t \in [-1, -1/2]$, $C_t$ is the same with $C_{1/2}$.  For $t \in [-3/2, -1)$, $C_t$ is the same as $C_{3/2}$, which is illustrated in (13).  Finally, the motion  picture for $C \cap B^2 \times [-2, -3/2]$ is the inverse of the motion picture for $C \cap B^2 \times [3/2, 2]$.  
}\end{example}

%
%

\section{From a curtain to a $3$-dimensional braid}\label{sect:curtain-3braid}

We have seen that a $3$-dimensional braid determines the braid monodromy, and the braid monodromy is described by a curtain.  In this section, we explain how to recover a $3$-dimensional braid, up to isomorphism, from its curtain description.  

Let $C$ be a curtain in $B^3 = B^2 \times [0,1]$ of degree $d$.   Let $\{ C_t \}_{t \in [0,1]}$ be the motion picture of $C$ and 
$E= \{ t_1, \dots, t_k\}$ the set of exceptional values.  For a regular value $t$, $C_t$ is a chart of degree $d$ describing a $2$-dimensional braid $S_t$ in $D^2 \times B^2$.  When a chart is deformed by an ambient  isotopy of $B^2$,  the corresponding $2$-dimensional braid is deformed by an ambient isotopy of $D^2 \times B^2$.  So we have a one-parameter family of equivalent $2$-dimensional braids.  For an exceptional value $t_j$ where a chart move occurs, the braid monodromy determined from the chart does not change, and the corresponding $2$-dimensional braids are isomorphic.  For an exceptional value $t_j$ where free edges are inserted, the corresponding $2$-dimensional braid is transformed into a $2$-dimensional braid by surgery along $1$-handles compatible with the structure of $2$-dimensional braids (cf. Chapter~20 of \cite{CK2012}.)  For an exceptional value $t_j$ where free edges are deleted, the inverse of surgery along $1$-handles (or equivalently, surgery along $2$-handles) occurs.  Considering the trace of these $2$-dimensional braids with surgeries along $1$-handles and $2$-handles, we have a motion picture of a $3$-dimensional braids in $D^2 \times B^3 = (D^2 \times B^2) \times [0,1]$. This is a desired $3$-dimensional braid. 

\begin{remark}\label{remark:closure}{\rm 
Assume that $S^3 = B^3_- \cup S^2 \times [0,1] \cup B^3_+$ where $\partial B^3_-$ is identified with $S^2 \times \{0\}$ and $\partial B^3_+$ is identified with $S^2 \times \{1\}$. 

Let $C$ be a curtain in $B^3 = B^2 \times [0,1]$ of degree $d$.   Let $\{ C_t \}_{t \in [0,1]}$ be the motion picture of $C$ and 
$E= \{ t_1, \dots, t_k\}$ the set of exceptional values. 
Assume that $B^2$ is contained in the $2$-sphere $S^2$ and then $B^3 = B^2 \times [0,1] \subset S^2 \times [0,1]$. 
For a regular value $t$, $C_t$ is regarded as a chart in $S^2$ describing a $2$-dimensional braid $S_t$ in $D^2 \times S^2$.  By the same argument as above, we have a one-parameter family of $2$-dimensional braids in $D^2 \times S^2$ in which surgeries along $1$-handles or $2$-handles may occur.  Since $S_0$ and $S_1$ are described by the empty charts $C_0$ and $C_1$ respectively, they are the trivial $2$-dimensional braid $X_d \times S^2$ in $D^2 \times S^2$.  The trace of this one-parameter family is a $3$-manifold embedded in $D^2 \times (S^2 \times [0,1])$.  Taking the union of this $3$-manifold with $X_d \times (B^3_- \cup B^3_+)$, we have a $3$-dimensional braid in $D^2 \times S^3$.  We call it a $3$-dimensional braid in $D^2 \times S^3$ described by the curtain $C$. Conversely, any $3$-dimensional braid in $D^2 \times S^3$ is equivalent to one obtained from a curtain this way.  
}\end{remark}

\begin{remark}{\rm 
Let $M$ be a $3$-dimensional braid in $D^2 \times S^3=  D^2 \times (B^3_- \cup S^2 \times [0,1] \cup B^3_+)$ and let $\Delta(M)$ be the set of branch values of the branched covering $M \to S^3$.  Deforming $M$ up to equivalence, we assume that $M$ is obtained as in Remark~\ref{remark:closure} from a curtain $C$ in $B^3 = B^2 \times [0,1]$.  
Suppose that $\Delta(M)$ is a braid form as in the proof of Theorem~\ref{thm:monodromy-curtain}. 
(Here we replace the values $t=-2, -1, 0, 1, 2$ in the proof with $t=0, 1/4, 1/2, 3/4, 1$, respectively.) 
Put $M_- = (M \cap D^2 \times (S^2 \times [0, 1/2])) \cup X_d \times B^3_-$ and 
$M_+ = (M \cap D^2 \times (S^2 \times [1/2, 1])) \cup X_d \times B^3_+$.  Then $M = M_- \cup M_+$.  
This splitting gives a Heegaard splitting of the $3$-manifold $M$.  
}\end{remark}

\begin{remark}{\rm 
In \cite{CK2012, CK*} the notion of $3$-dimensional braids are extended to immersed $3$-manifolds.  For immersed $3$-dimensional braids, braid monodromies and curtain descriptions can be also considered.  
Further detailed research on these descriptions is expected.   
}\end{remark}

\begin{remark}{\rm 
In this paper we assume that the branch set is a link in $B^3$.  For a definition of a $4$-dimensional braid, it seems to be natural to consider the branch set is (i) an embeddded surface or (ii) an immersed surface with transverse double points (cf. \cite{IP2002, Pi1995}). An example of a chart description for a $4$-dimensional braid is given in \cite{CK*}. 
}\end{remark}


\begin{thebibliography}{99}



%

\bibitem{BE1979} I. Berstein and A. L. Edmonds,  {\it On the construction of branched coverings of low-dimensional manifolds}, Trans. Amer. Math. Soc. {\bf 247}  (1979), 87--124. 
%

\bibitem{BE1984} I. Berstein and A. L. Edmonds, {\it On the classification of generic branched coverings of surfaces}, Illinois J. Math. {\bf 28} (1984), no. 1, 64--82.  

\bibitem{CK2012}
J. S. Carter and S. Kamada, {\it Braids and branched coverings of dimension three}, 
in the proceedings of \lq\lq Intelligence of low dimensional topology\rq\rq, (1012), RIMS Kokyuroku {\bf 1817}, pp. 64--82,
arXiv:1206.4744v1. 

\bibitem{CK*}
J. S. Carter and S. Kamada, {\it How to fold a manifold}, 
arXiv:1301.4259v1. 


\bibitem{CKS2004}
J. S. Carter, S. Kamada and M.  Saito, {\it Surfaces in $4$-space}, Encyclopaedia of Mathematical Sciences, 142. Low-Dimensional Topology, III. Springer-Verlag, Berlin, 2004.  

\bibitem{CS1998} J. S. Carter and M. Saito, {\it Knotted surfaces and their diagrams},  Math. Surveys Monogr. {\bf 55},  American Mathematical Society, Providence, RI, 1998.  





\bibitem{Hi1976}
H. M. Hilden, {\it Three-fold branched coverings of $S^3$}, Amer. J. Math. {\bf 98} (1976), 989--997.  
%


\bibitem{Hud1969}
J. F. P. Hudson, {\it Piecewise linear topology}, Benjamin, New York, 1969. 

%
%

\bibitem{IP2002}
M. Iori and R. Piergallini, {\it $4$-manifolds as covers of the $4$-sphere branched over non-singular surfaces}, Geom. Topol. {\bf 6} (2002), 393--401. 
%


\bibitem{Kam1992}
S. Kamada, 
{\it Surfaces in $R^4$ of braid index three are ribbon}, 
J. Knot Theory Ramifications {\bf 1} (1992) 137--160.

\bibitem{Kam1996} S. Kamada,  {\it An observation of surface braids via chart description},  J. Knot Theory Ramifications {\bf 4} (1996), 517--529. 

\bibitem{Kam2002}
S. Kamada, {\it Braid and knot theory in dimension four}, Math. Surveys Monogr. {\bf 95}, 
 American Mathematical Society, Providence, RI, 2002.







\bibitem{Mo1976}
J. M. Montesinos, {\it Three-manifolds as 3-fold branched covers of $S^3$}, Quart. J. Math. Oxford Ser. (2) {\bf 27} (1976), 85--94. 
%

%

\bibitem{Pi1995}
R. Piergallini, {\it  Four-manifolds as $4$-fold branched covers of $S^4$},  Topology {\bf 34} (1995), no. 3, 497--508. 
%

\bibitem{RS1972}
C. P. Rourke and B. J. Sanderson, {\it Introduction to piecewise-linear topology}, Ergebnisse der Mathematik und ihrer Grenzgebiete, Band 69. Springer-Verlag, New York-Heidelberg, 1972.  

\bibitem{Ta2007}
K. Tanaka, {\it A note on CI-moves}, in Intelligence of Low dimensional topology 2006, pp. 307--314, Ser. Knots Everything, {\bf 40}, World Scientific Publishing Co., 2007. 


\end{thebibliography}
\end{document}